\documentclass[12pt]{elsarticle}
\usepackage{amsmath}
\usepackage{amssymb}
\usepackage[utf8]{inputenc}
\usepackage[T1]{fontenc}
\usepackage[english]{babel}
\usepackage{color}
\usepackage{booktabs}
\usepackage{tabularx}
\usepackage{array}
\usepackage{rotating}
\journal{Physica D}

\begin{document}

\begin{frontmatter}

\title{From Complex Regularization to Polynomial Flow: Capture and Escape in Rocking Dynamics}

\author[inst1]{Francisco Prieto-Castrillo}
\affiliation[inst1]{organization={Departamento de Matemáticas, Universidad de Oviedo},%Department and Organization
            addressline={Avda Calvo Sotelo s/n}, 
            city={Oviedo},
            postcode={33007}, 
            country={Spain}}

\author[inst1]{Set Pérez-González}
\author[inst1]{Jesús S. Pérez del Río}

\begin{abstract}
We develop a planar polynomial system based on a complex regularization of free rocking and analyze the equilibrium states and phase portraits via an effective radial potential. The conflict between the centrifugal effect, which is introduced by means of regularization, and gravity yields a potential well and an escape barrier for $0<\beta<27/256$; both vanish at a saddle--node bifurcation point. The exact energy conservation law allows us to separate dissipative relaxation, conservative periodic motions along the radial direction, and energy-creating trajectories. The node-focus boundary of finite equilibria is different from the degeneracy condition for compactification at infinity. As the last step of analysis, we reconstruct the complex phase portrait and discover that a stable radial equilibrium corresponds to a relative equilibrium in the regularized dynamics, whereas the reconstruction of an impact dynamics requires taking a singular limit.
\end{abstract}


\begin{highlights}
\item Polynomial coordinates describe regularized radial rocking dynamics.
\item Energy balance distinguishes radial capture, conservative motion, and escape.
\item Phase reconstruction clarifies the scope of the physical impact limit.
\end{highlights}

\end{frontmatter}

%%%%%%%%%%%%%%%%%%%%%%%%%%%%
%%%%%%%%%%%%%%%%%%%%%%%%%%%%
\section{Rocking: a ubiquitous and intricate dynamic behavior}
\label{sec:intro}
A passenger on a train moving in a particular direction and suddenly applying brakes finds himself having to shift his body weight, placing his feet ahead and turning, so that he maintains his balance. As a result, his center of gravity moves slightly backward, making him turn on the other foot. Such an action is a controlled oscillation that serves to stabilize oneself. This is one of the instances of rocking motion (RM).
\begin{figure}[h!]
\begin{center}
\includegraphics[width=\linewidth]{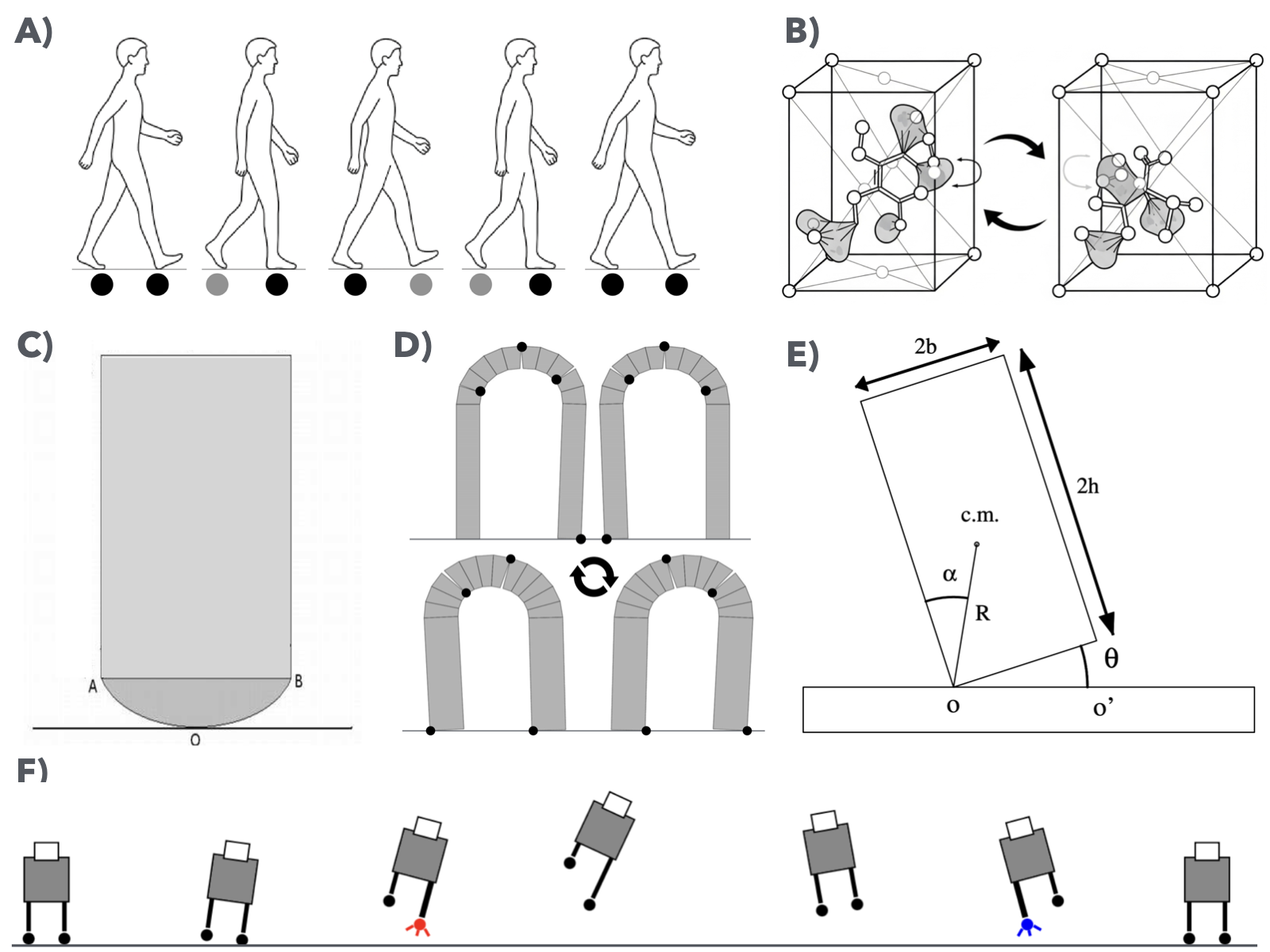} 
\caption{Some examples of rocking-type movements: (A) a walk in which the rotational axes switch dynamically to create a successful walk; (B) a change of the molecule from one conformational state to another, both of which have a specific effective rotation axis in the crystal lattice; (C) a curved-bottomed block for which the instantaneous rotational axis changes continuously on the line of contact; (D) two types of collapse of four-hinged masonry arches, in which some of the structural elements have sub-rocking; (E) a rigid block rocking back and forth on its corners; (F) a robot rocking itself forward over one of its rotational points using a rocking mechanism. In all these systems, the feature they share is that they oscillate from one stable position to another through rotating around different axes. These rotations can be either continuous (A and C) or discontinuous (B, D, E and F).}
\label{fig:draft_rm_universal}
\end{center}
\end{figure}
Indeed, this approach is synergistically combined with others for achieving propulsion in human and robotic locomotion \cite{Walking2006,Walking2019, Raibert1986}. During each step, a dynamic change takes place from one axis of instantaneous rotation to another, such that the leading leg becomes the new pivot leg and the trailing leg gets detached, starting a series of controlled oscillations for facilitating the forward movement. Such an alternation of axes of rotations between legs is a simple but necessary example of rocking during the gait cycle \cite{Walking2006,Walking2019}. Moreover, besides its application in locomotion, this very mechanism defines the functioning of some elements of ancient constructions in the presence of earthquake loads \cite{Makris2013}. Rocking motion is seen in macromolecules, helping their internal dynamics within the crystalline form \cite{Proteins2015,Proteins2023,MolecularRotor}, as well as in meta-constructions with highly desired mechanical properties \cite{MetaInterfaces} (see Fig.~\ref{fig:draft_rm_universal}).
\par 
The similarity among these diverse dynamic systems is that they oscillate about an equilibrium position through rotations along one or more axes. An important distinction is the type of transition between these rotational axes, which could be gradual or abrupt. If the transitions between rotational axes are sufficiently gradual, it is often possible to treat the dynamics of the system as a variation of the harmonic oscillator, particularly when there are small deviations (see Fig.\ref{fig:draft_rm_universal} C). On the other hand, if the transitions are abrupt, it might seem natural to use piecewise equations, whereby each equation corresponds to the dynamics of rotation about a particular axis. Moreover, these abrupt transitions could be accompanied by impulsive actions that either add or subtract energy from the system.
%RM complejo
\par 
However, in the second case, when there is a combination of non-smooth dynamics with impulsive loads, the analysis becomes more complicated. Indeed, in some cases of this kind of dynamics, very complicated effects arise, such as chaos \cite{Hogan,Rega2006} and chattering \cite{Chattering2001}. Such dynamics require a very sophisticated mathematical approach, which is one of the reasons for the development of scientific research in this area. As an example of such scientific investigation, we can mention successful integration of rocking dynamics analysis and artificial intelligence \cite{AI}.
\par 
%Meccanica
The use of the rigid-block model -Fig.\ref{fig:draft_rm_universal} E- is an example of an impulsive system of a sudden and dissipative type, which has been transformed to a central force field problem in the complex plane \cite{Meccanica}. In \cite{Meccanica} the authors found an explicit parameter that regulated the smoothness of the dynamics, and the RM limit turned out to be a special case of the central force field problem. Moreover, impulsive forces that were initially added ad hoc were nicely incorporated into the new approach as generalized forces. The main advantage of such reformulation of the problem was the possibility to describe the dynamics with Ordinary Differential Equations (ODEs) instead of Piecewise Smooth (PWS) systems at the expense of introducing singularities in the equations.
\par 
%This paper
Following up on this advance, we generalize the regularization scheme of \cite{Meccanica} to take into account the dynamics of the singularities properly, which gives rise to polynomial dynamic equations and thereby renders the system susceptible to planar treatment. In addition, making use of the broad applicability of the theory proposed by \cite{Meccanica}, our findings are more general than the extreme situation studied in their work.
In fact, we determine the criteria for radial confinement and radial escape, classify monotonic and oscillatory relaxation and delineate the interpretation of the generalized parameter plane.
 %Organizacion del paper
\par 
The paper is structured as follows.
First we provide a reminder about the regularization of \cite{Meccanica}.
Then we consider polynomials of coordinates and examine finite equilibria, periodic solutions, and local behavior at infinity through the Poincaré compactification method.
The last two sections interpret the previous findings in the context of radial mechanics and provide conclusions for their implications and limitations to  rocking systems.
%%%%%%%%%%%%%%%%%%%%%%%%%%%%
\section{Complex regularization}
\label{sec:compregu}
%%%%%%%%%%%%%%%%%%%%%%%%%%%%
As in \cite{Meccanica}, we will use the abrupt and impulsive-dissipative model of the rigid block as our starting point. We summarize here the formulation and the main results from \cite{Meccanica} for reference.
\par
\begin{figure}[h]
 \begin{center}
 \includegraphics[scale=0.4]{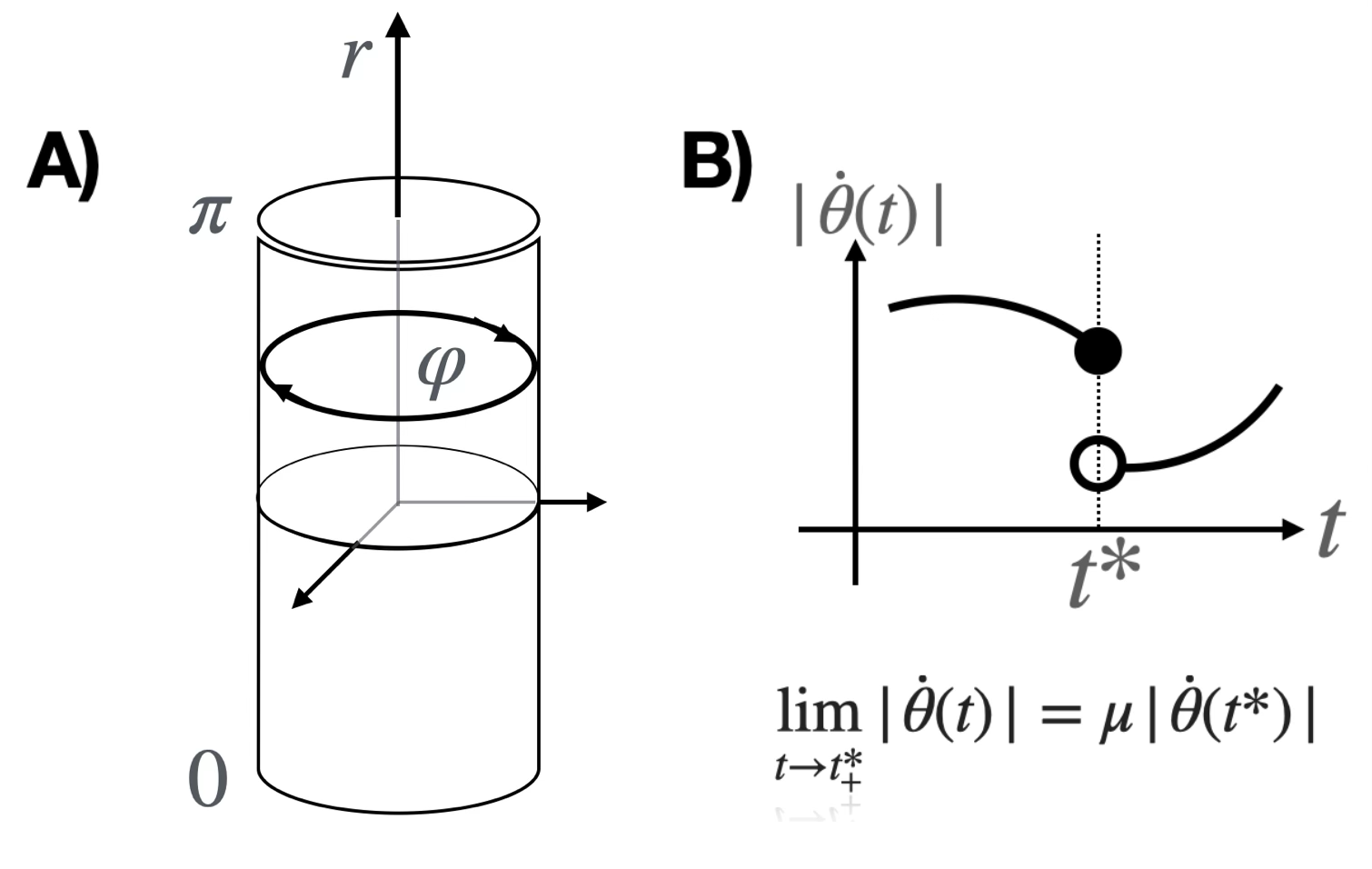}	
 \end{center}
\caption{(A) Regularization of the rocking angle; (B) behavior of the angular velocity at impact.}
\label{fig:draft_regularization}
 \end{figure}

 
Let $\theta$ denote the rocking angle of the block depicted in Fig.~\ref{fig:draft_rm_universal}E. 
The configuration space is $\mathcal{Q}_1 = \mathbb{S}^1$, representing the rocking angle $\theta$ as an angular coordinate. The tangent space is $T\mathcal{Q}_1 = \mathbb{S}^1 \times \mathbb{R}$, with elements $(\theta, \dot{\theta})$. The phase space is $T^*\mathcal{Q}_1 = \mathbb{S}^1 \times \mathbb{R}$, with elements $(\theta, p_\theta)$.
The Lagrangian  $L: T\mathcal{Q}_1 \to \mathbb{R}$  of this system is $L(\theta,\dot{\theta})= (I/2)\dot{\theta}^2 -MgR\cos(\theta\mp\alpha)$, where the upper and lower signs refer to the domains $\theta>0$ and $\theta<0$, respectively; $I$ is the moment of inertia about the center(s) of rotation, $M$ is the mass of the object, $R$ is the distance from the center of rotation to the center of mass, $g$ is the gravitational acceleration, and $\alpha=\arctan(b/h)$ is the so-called \textit{slenderness} angle. The Euler--Lagrange equation for $\theta$, together with the restitution relation for the angular velocity $\dot{\theta}$ (see Fig.~\ref{fig:draft_regularization}B), yields
\begin{equation}\label{eq:PFRM}
\begin{array}{r}
I\ddot{\theta}+MgR\sin(\alpha\operatorname{sgn}\theta-\theta)=0,\quad\theta\neq0\\
|\dot\theta(t^{*+})|=\mu|\dot\theta(t^{*-})|,\quad\mu\in[0,1]
\end{array}
\end{equation}
Here $t^{*}$ is one of the impact times (in general, there is a sequence of them). We refer to these equations as the piecewise free rocking motion (PFRM) equations.\\
Our first regularization consists in extending the rocking angle from $\mathbb{R}$ to $\mathbb{C}$ through the definition: $\theta = re^{i\psi}$, where $r>0$ and $\psi\in\mathbb S^1$. Away from the origin the polar configuration space is $(0,\infty)\times\mathbb S^1$ (see Fig.~\ref{fig:draft_regularization} A).
 The configuration space is $\mathcal{Q}_2 = (0,\infty)\times\mathbb S^1$. The tangent space is $T\mathcal{Q}_2 = \mathcal Q_2 \times \mathbb{R}^2$, with elements $(r, \psi, \dot{r}, \dot{\psi})$. The phase space is $T^*\mathcal{Q}_2 = \mathcal Q_2 \times \mathbb{R}^2$, with elements $(r, \psi, p_r, p_\psi)$.
 Following \cite{Meccanica}, the simplest real kinetic energy is $T = \frac{I}{2}\dot\theta\dot{\theta}^*$. Together with the gravitational potential $U(r)=MgR\cos(r-\alpha)$, it leads to the Lagrangian $L=T-U: T\mathcal{Q}_2 \to \mathbb{R}$ in Eq.~\eqref{eq:lagrangian}.
 \begin{equation}\label{eq:lagrangian}
 	L(r,\dot{r},\psi,\dot\psi) = \frac{I}{2}(\dot r^2+r^2\dot\psi^2)-MgR\cos(r -\alpha)
 \end{equation}
The generalized momentum $p_\psi = \frac{\partial L}{\partial \dot{\psi}} = Ir^2\dot{\psi}$ is conserved due to $\psi$ being a cyclic coordinate.
The conservative Hamiltonian $H: T^*\mathcal{Q}_2 \to \mathbb{R}$ is obtained via Legendre transformation:
\begin{equation}\label{eq:hamiltonian}
H(r, p_r, p_\psi) = \frac{p_r^2}{2I} + \frac{p_\psi^2}{2Ir^2} + MgR\cos(r - \alpha)
\end{equation}
 This describes motion in an effective central potential, obtained by adding an infinite repulsive barrier at $r=0$ to $U(r)$.\\
 To implement the impact forces, we build on the results reported in \cite{Meccanica,Dirac} that modeled the impact as the simplest coupling between the two velocities involved $\dot r$ and $\dot \psi$. 
 In particular, it was shown in \cite{Meccanica} that the phase $\psi(t)$ exhibits the following limiting behavior:
 \begin{equation}
 	\lim_{p_\psi\to0}\psi(t)=\pi\sum_{j}\Theta(t-t_j^*)
 	\label{eq:phase_dynamics}
 \end{equation}
where $\Theta$ is the Heaviside function and $t_j^*$ are impact times. The radial generalized force couples radial and phase velocities. With the dimensionless coupling convention used here,
\begin{equation}
Q_r^d=I\sigma\dot r\dot\psi,\qquad
\sigma=\frac{2}{\pi}\ln\mu,\qquad \mu>0,
\end{equation}
where the restitution calibration is inherited from \cite{Meccanica}. Its power is $Q_r^d\dot r=I\sigma\dot r^2\dot\psi$, nonpositive for $\sigma\leq0$ and $p_\psi>0$. The radial equation follows from
\begin{equation}
I(\ddot r-r\dot\psi^2)=Q_r^U+Q_r^d.
\end{equation}
The centrifugal term must be retained when analyzing the impact layer: the logarithmic restitution relation is a singular-limit calibration, not the result of dropping $r\dot\psi^2$ in this equation. In what follows we analyze the smooth system with that prescribed calibration.
\par
 The conservative generalized force on $r$ is $Q_r^U=-\frac{\partial U}{\partial r}=MgR\sin(r-\alpha)$.\\
 Combining $Q_r^d$ and $Q_r^U$, the radial equation of motion becomes
 \begin{equation}
 	\ddot r=r\dot\psi^2+k\sin(r-\alpha)+(2\ln\mu/\pi)\dot r\dot\psi
 \end{equation}
 where we have introduced the quantity $k\equiv MgR/I$. The equation for $\psi$ is simply the conservation law $\dot p_\psi = 0$.\\
 This equation is now written in normalized form by defining $\rho\equiv r/\alpha$ and the dimensionless time $\tau\equiv t\sqrt{k}$:
 \begin{equation}
 	\rho''=\rho(\psi')^2+(1/\alpha)\sin(\alpha \rho-\alpha)+(2\ln\mu/\pi)\rho'\psi'
 \end{equation}
 with the notation $f'\equiv df/d\tau$.\\
 Substituting $\psi'=p/\rho^2$ and rearranging gives
 \begin{equation}
 	\rho''=\frac{p^2}{\rho^3}-(1/\alpha)\sin(\alpha -\rho\alpha)+\biggl(\frac{2p\ln \mu }{\pi }\biggr)\frac{\rho'}{\rho^2}
 \end{equation}
where we have defined the dimensionless parameter $p\equiv p_\psi\alpha^{-2}(MgRI)^{-1/2}$.
Linearizing the sine for $|\alpha(1-\rho)|\ll1$, so that $\sin(\alpha-\alpha\rho)\approx\alpha(1-\rho)$, produces the small-angle approximation:
\begin{equation}\label{eq:Rockingfinal}
	\rho''=p\sigma \frac{\rho'}{\rho^2}+\biggl(\rho-1+\frac{p^2}{\rho^3}\biggr)
\end{equation} 
 For small $p$, this regularization reduces to the PFRM equations~\eqref{eq:PFRM}, in excellent agreement with theory \cite{Meccanica} and experiment \cite{Pena}.
%%%%%%%%%%%%%%%%
\section{Desingularization through Poincaré compactification: a polynomial model}
%%%%%%%%%%%%%%%%%
Instead of \eqref{eq:Rockingfinal}, in this section we consider the following generalized system
\begin{equation}\label{sys:Rockinggeneral}
\begin{cases}
\rho'=\omega\\
\omega'=\gamma \frac{\omega}{\rho^2}+\biggl(\rho-1+\frac{\beta}{\rho^3}\biggr)
\end{cases}
\end{equation}
where $\gamma$ and $\beta$ are arbitrary real numbers. Note that this generalization includes cases outside physical rocking motion \eqref{eq:Rockingfinal}, which corresponds to $\gamma=p \sigma$ and $\beta=p^2$.

To gain a comprehensive understanding of the behavior of the planar system, we analyze the solutions that tend toward or originate from infinity. This requires a compactification of the plane. Several compactification methods exist in the literature, with the most common being those that represent infinity either as a point or as a circle. In this work, we adopt the Poincaré compactification, which maps a planar vector field $X$ onto a vector field defined on the sphere. 

Moreover, to study the qualitative behavior of system \eqref{sys:Rockinggeneral}, we would like firstly to remove its singularities. The Poincaré compactification can be used as a change of variables that aids desingularization; see \cite{GarciaPerezChavelaSusin2006}.

\begin{figure}
\centering
\includegraphics[width=0.56\linewidth]{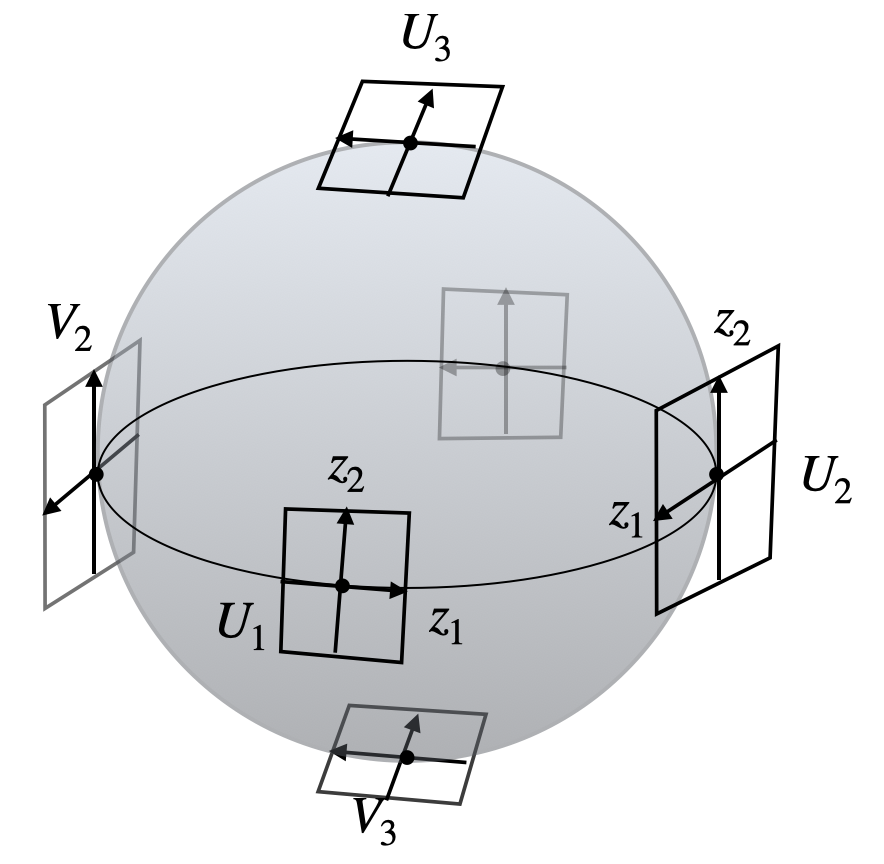}
\caption{Local charts on the sphere in the Poincaré compactification.}\label{fig:cartasP}
\end{figure}

This procedure uses six local charts, centered at the poles and on the equator, which together cover the whole sphere (see Figure~\ref{fig:cartasP}).

\begin{figure}
\centering
\includegraphics[width=0.6\linewidth]{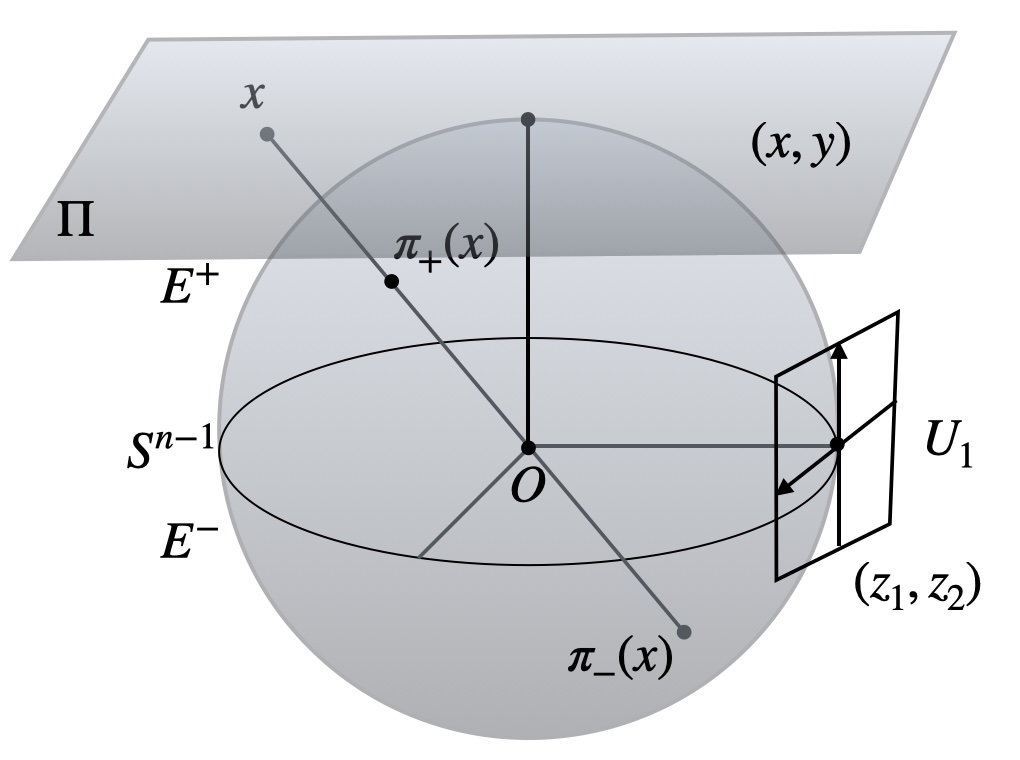}
\caption{Projection onto the sphere in the Poincaré compactification.}
\end{figure}

Consider the vector field $X$ defined by
\begin{equation}\label{sys:VFPoincare}
\begin{cases}
u_1' = P(u_1,u_2)\\
u_2' = Q(u_1,u_2).
\end{cases}
\end{equation}
on the real plane embedded in \(\mathbb{R}^3\) as the tangent plane to the unit sphere \(S^2\) at the north pole. Accordingly, points on the plane can be represented as \((u_1,u_2, 1)\). Each such point can be projected onto the upper hemisphere by connecting it with the origin through a straight line. If \((z_1, z_2, z_3)\) are the coordinates of the corresponding point on the sphere, \(S^2\), then this projection satisfies \( u_1= z_1 / z_3\) and \(u_2 = z_2 / z_3\). Conversely, the coordinates on the sphere can be expressed as
\[
z_1 = \frac{u_1}{\Delta}, \quad z_2 = \frac{u_2}{\Delta}, \quad z_3 = \frac{1}{\Delta},
\]
where \(\Delta = \sqrt{1 + u_1^2 + u_2^2}\). The projection onto the southern hemisphere is analogous.

To analyze the critical points at infinity, it is necessary to employ local charts. In the chart \( U_1=\{z\in S^2/z_1>0\} \), the change of variables is given by \( (u_1, u_2) = \left( \frac{1}{z_2}, \frac{z_1}{z_2} \right) \). Moreover, in the chart \( U_2=\{z\in S^2/z_2>0\}  \), the change of variables is given by \( (u_1, u_2) = \left( \frac{z_1}{z_2}, \frac{1}{z_2} \right) \). Finally, in the chart \( U_3=\{z\in S^2/z_3>0\}  \), the change of variables is given by \( (u_1, u_2) = \left( z_1, z_2 \right) \). Straightforward computations show that, in these charts, $X$ is given by the systems
\begin{equation}\label{sys:cartaU1}
\begin{cases}
z_1' = z_2^n\left(-z_1 P\left(\frac{1}{z_2},\frac{z_1}{z_2}\right)+Q\left(\frac{1}{z_2},\frac{z_1}{z_2}\right)\right) \\
z_2' = -z_2^{n+1}P\left(\frac{1}{z_2},\frac{z_1}{z_2}\right).
\end{cases}
\end{equation}

\begin{equation}\label{sys:cartaU2}
\begin{cases}
z_1' = z_2^n\left(P\left(\frac{z_1}{z_2},\frac{1}{z_2}\right)-z_1Q\left(\frac{z_1}{z_2},\frac{1}{z_2}\right)\right),\\
z_2' = -z_2^{n+1}Q\left(\frac{z_1}{z_2},\frac{1}{z_2}\right).
\end{cases}
\end{equation}

\begin{equation}\label{sys:cartaU3}
\begin{cases}
z_1' = P(z_1,z_2)\\
z_2' = Q(z_1,z_2).
\end{cases}
\end{equation}
 For the nonpolynomial field \eqref{sys:Rockinggeneral} we take $n=1$: this choice involves no time rescaling and already yields a polynomial system.

Under the Poincar\'e compactification, the dynamical behavior of $X$ near infinity corresponds to the behavior in a neighborhood of the equator $\mathbb S^1$ of $S^2$. In the local chart $U_1$, the singular points with $z_2 = 0$ are the infinite singular points of this chart. In the local chart $U_2$, the only remaining candidate is the origin, $z_1 = z_2 = 0$. For a more thorough account of this compactification, we refer the reader to \cite{DumortierLlibreArtes2006}.

Applying the change of variables of chart $U_1$, Eq.~\eqref{sys:cartaU1}, to the singular system \eqref{sys:Rockinggeneral}, we obtain
\begin{equation}\label{sys:nuestracartaU1}
\begin{cases}
x' = -x^2+\gamma x y^2+1-y+\beta y^4\\
y' = -xy,
\end{cases}
\end{equation}
The explicit change, with no time rescaling, is
\begin{equation}\label{eq:physicalcoordinates}
x=\frac{\omega}{\rho},\qquad y=\frac{1}{\rho},\qquad
\rho=\frac{1}{y},\qquad \omega=\frac{x}{y}.
\end{equation}
Thus the polynomial system is equivalent to the smooth radial model on $y>0$. Its extension to $y\leq0$ and its compactification provide mathematical boundary information; their interpretation as impacting motion requires the reconstruction discussed in Section~\ref{sec:physical}.

\section{Qualitative behavior}

In this section we study the qualitative behavior of system \eqref{sys:nuestracartaU1}. We analyze the singular points and their bifurcations (Subsection~\ref{sec:singular}), the nonexistence of limit cycles (Subsection~\ref{sec:nocycles}) and the local dynamics at infinity (Subsection~\ref{sec:infinity}). These results are then summarized in the parameter diagram and illustrated by computed phase portraits (Subsections~\ref{sec:paramdiagram} and~\ref{sec:compactified}).
\subsection{Singular points}\label{sec:singular}
We begin by determining the singular points of system \eqref{sys:nuestracartaU1}. The second equation implies that every singular point lies on one of the lines $x=0$ or $y=0$. On $y=0$ there are exactly two singular points, $P_1=(-1,0)$ and $P_2=(1,0)$. On $x=0$, the singular points are given by the solutions of the equation
\begin{equation}\label{eq:despejabeta}
1-y+\beta y^4=0.
\end{equation}
Each solution $\tilde{y}$ of this equation corresponds to an intersection of the graphs of $1-y$ and $-\beta y^4$. The number of real solutions depends on $\beta$, as detailed below.

For particular values $\beta^*$ and $y^*$, which must be determined explicitly, the two graphs are tangent (double contact). To compute these values, we solve the following system for the unknowns $y$ and $\beta$:
\begin{equation}
\left\{\begin{array}{l}
1-y+\beta y^4=0\\
\frac{d}{dy}\left(-\beta y^4\right)=-1
\end{array}  \right..
\end{equation}
Its solution is $\beta^*=\frac{3^3}{4^4}=\frac{27}{256}$ and $y^*=\frac{4}{3}$.

Hence, depending on the value of $\beta$:
\begin{enumerate}
\item if $\beta>\beta^*$, there are no singular points other than $P_1$ and $P_2$;
\item if $\beta=\beta^*$, there is one additional real equilibrium, $P_3=(0,y^*)$;
\item if $\beta\in(0,\beta^*)$, there are two additional equilibria, $P_3=(0,\tilde{y})$ and $P_4=(0,\hat{y})$, with $\hat{y}>y^*>\tilde{y}>1$;
\item if $\beta=0$, there is one additional real equilibrium, $P_3=(0,1)$;
\item if $\beta<0$, there are two additional equilibria, $P_3=(0,\tilde{y})$ with $\tilde{y}\in(0,1)$, and $P_4=(0,\hat{y})$ with $\hat{y}<0$.
\end{enumerate}

To continue the qualitative study of system \eqref{sys:nuestracartaU1}, we consider its Jacobian matrix,
\begin{equation}
J(x,y)=\begin{pmatrix}
-2x+\gamma y^2&2\gamma xy-1+4\beta y^3\\
-y&-x
\end{pmatrix}
\end{equation}
At $P_1$ and $P_2$ the eigenvalues of $J$ are $\lambda_1=-2x$ and $\lambda_2=-x$; consequently, $P_1=(-1,0)$ is a repelling node and $P_2=(1,0)$ is an attracting node. For the equilibria $P_3$ and $P_4$, which are of the form $(0,y)$, Eq.~\eqref{eq:despejabeta} can be used to simplify the eigenvalues of the Jacobian matrix, which become
\begin{equation}
\lambda_{1,2}=\frac{\gamma {y}^2\pm \sqrt{\gamma^2 {y}^4-4(3{y}-4)}}{2}.
\end{equation}
Thus, for $y<y^*= \frac{4}{3}$ both eigenvalues are real and of opposite signs, so all equilibria located below the line $y=\frac{4}{3}$ are saddles.

If $y=y^*$, then $(0,y^*)$ is an equilibrium when $\beta=\beta^*$. In this case,
\begin{equation}
J(0,y^*)=\begin{pmatrix}
 \frac{16}{9}\gamma&0\\
-\frac{4}{3}&0
\end{pmatrix}.
\end{equation}
Therefore, following the procedure described in \cite{DumortierLlibreArtes2006}, if $\gamma=0$ the equilibrium is nilpotent and of cusp type (see case 5 in Figure~\ref{fig:phasecomputed}). Moreover, if $\gamma \neq 0$, the point is semi-hyperbolic: it is an attracting (respectively, repelling) saddle--node when $\gamma<0$ (respectively, $\gamma >0$); see cases 2 and 8 in Figure~\ref{fig:phasecomputed}.

If $y>y^*= \frac{4}{3}$, the type of the equilibrium depends on $\gamma$. If $\gamma^2 {y}^4-4(3{y}-4)=0$, then $\lambda_1=\lambda_2\in\mathbb{R}$, with the same sign as $\gamma$. In this case the equilibrium is an improper node: since $J_{21}=-y\neq0$, the matrix is not a multiple of the identity. If $\gamma^2 {y}^4-4(3{y}-4)>0$, the eigenvalues are real, distinct and of the same sign, again that of $\gamma$, and the equilibrium is a node. Finally, if $\gamma^2 {y}^4-4(3{y}-4)<0$, the eigenvalues are complex conjugate and the equilibrium is a focus or a center. The center corresponds to $\gamma=0$ (the first integral of Subsection~\ref{sec:nocycles} shows that it is a genuine nonlinear center, not only a linear one); otherwise, the focus is attracting or repelling according to the sign of $\gamma$.

For the bifurcation diagram in the parameters $(\gamma,\beta)$, we need the resultant of Eq.~\eqref{eq:despejabeta} and $\gamma^2 {y}^4-4(3{y}-4)=0$, regarded as polynomials in $y$:
\begin{equation}
\begin{aligned}
\mathrm{Res}={}&65536\beta^4-256(64\gamma^2+27)\beta^3\\
&+192(8\gamma^4+9\gamma^2)\beta^2\\
&-16(4\gamma^6+9\gamma^4)\beta+4\gamma^6+\gamma^8.
\end{aligned}
\end{equation}
%\begin{align}
%Res&=\gamma^8 - 64\gamma^6\beta + 4\gamma^6 + 1536\gamma^4\beta^2 - 144\gamma^4\beta - 16384\gamma^2\beta^3 + 1728\gamma^2\beta^2 + 65536\beta^4 - 6912\beta^3\\
%&=\gamma^8 + (- 64\beta + 4)\gamma^6 + (1536\beta - 144)\gamma^4\beta + (- 16384\beta + 1728)\gamma^2\beta^2 + (65536\beta - 6912)\beta^3\\
%&=\gamma^8 + 4(- 16\beta + 1)\gamma^6 + 48(32\beta - 3)\gamma^4\beta + 64(-256\beta + 27)\gamma^2\beta^2 + 256(256\beta - 27)\beta^3\\ 
%&=65536\beta^4 + (- 16384\gamma^2 - 6912)\beta^3 + (1536\gamma^4 + 1728\gamma^2)\beta^2 + (- 64\gamma^6 - 144\gamma^4)\beta + \gamma^8 + 4\gamma^6\\
%&=65536\beta^4 -256 (64\gamma^2 +27)\beta^3 + 192(8\gamma^4 + 9\gamma^2)\beta^2 -16(4\gamma^6 + 9\gamma^4)\beta + \gamma^8 + 4\gamma^6
%\end{align}

\subsection{Nonexistence of limit cycles}\label{sec:nocycles}
For $\gamma\neq0$, the ordinary divergence of \eqref{sys:nuestracartaU1} is $-3x+\gamma y^2$, which has no fixed sign. Instead, the Dulac multiplier $B=y^{-3}$ gives
\begin{equation}
\frac{\partial(BP)}{\partial x}+\frac{\partial(BQ)}{\partial y}
=\frac{\gamma}{y}.
\end{equation}
This expression has a strictly constant sign on each simply connected half-plane $y>0$ and $y<0$, excluding periodic orbits there. A periodic orbit cannot cross the invariant line $y=0$, and the scalar flow on that line cannot be periodic.

For $\gamma=0$, the first integral on $y\neq0$ is
\begin{equation}
\mathcal H(x,y)=\frac{x^2+2y-1+\beta y^4}{y^2}.
\end{equation}
Any regular periodic level is surrounded by nearby periodic levels and is therefore not an isolated limit cycle. For $0<\beta<\beta^*$, these levels form the period annulus around the center. Thus there are no limit cycles for any real $\beta,\gamma$, although conservative periodic radial motions do exist. An equivalent physical-domain proof follows from the energy balance in Section~\ref{sec:physical}.

\subsection{Local dynamics at polynomial infinity}\label{sec:infinity}
Infinity of the $(x,y)$ plane must be distinguished from $\rho\to\infty$, which is represented by $y\to0$. In particular, $y\to+\infty$ corresponds to $\rho\to0^+$, whereas $|x|\to\infty$ can also encode unbounded radial velocity. The inverse map \eqref{eq:physicalcoordinates}, rather than the compactified radius alone, determines the physical meaning of a boundary approach.

For $\beta\neq0$ the polynomial degree is four. In the chart $(x,y)=(1/v,u/v)$ of Poincarè's compactification, the transformed vector field gives
\begin{equation}\label{sys:nuestracartaU1deU1}
\begin{cases}
\dot u=u(-v^4+uv^3-\gamma u^2v-\beta u^4),\\
\dot v=v(v^2-v^4+uv^3-\gamma u^2v-\beta u^4).
\end{cases}
\end{equation}
Dots in this subsection refer to desingularized time, not physical time. The only equilibrium on $v=0$ in this chart is $u=0$, since the highest-degree angular polynomial is $-\beta u^5$. The complementary vertical direction has angular polynomial $\beta$ and contributes no additional equilibrium.

The weighted substitution $u=w$, $v=w^2z$, followed by removal of the common factor $w^4$, yields
\begin{equation}\label{eq:desingularizadoblowup}
\begin{cases}
\dot w=w(-\beta-\gamma z+w^3z^3-w^4z^4),\\
\dot z=z(z^2+\beta+\gamma z-w^3z^3+w^4z^4).
\end{cases}
\end{equation}
On the exceptional divisor $w=0$, the equilibria satisfy $z(z^2+\gamma z+\beta)=0$. At $A=(0,0)$ the eigenvalues are $-\beta,\beta$, so $A$ is a saddle. Additional real roots exist when $\Delta=\gamma^2-4\beta\geq0$:
\begin{equation}
z_\pm=\frac{-\gamma\pm\sqrt{\Delta}}{2},\qquad
\lambda_w=z_\pm^2,\qquad \lambda_z=\pm z_\pm\sqrt{\Delta}.
\label{eq:autovaloresBlowup}
\end{equation}
For $\beta>0$, $\Delta>0$, there is one source and one saddle in this desingularized chart, with their assignment to $z_+$ and $z_-$ depending on the sign of $\gamma$. For $\beta<0$, $z_+>0>z_-$ and both additional equilibria are sources. At $\Delta=0$ the tangential eigenvalue vanishes, so the hyperbolic classification does not apply. These statements describe this blow-up chart. Time orientation and the other directional charts must be accounted for asserting a complete global classification.

Indeed, in these coordinates
\begin{equation}
\rho=wz,\qquad \omega=\frac{1}{w},\qquad z=\rho\omega.
\end{equation}
Approach to $w=0$ at a nonzero finite $z$ combines $\rho\to0$ with unbounded radial velocity. Moreover, the factor $v^3$ reverses time orientation for $v<0$. Consequently, a source in the desingularized vector field does not establish physical uplift from an impact. These branches would be singular asymptotic regimes of the extended equations, not finite-velocity impact laws.

The condition $\gamma^2=4\beta$ concerns the exceptional divisor. It is different from the finite equilibrium node--focus boundary derived below. In addition, for $\beta=0$, the polynomial degree drops to three if $\gamma\neq0$ and to two if $\gamma=0$. The degree four boundary analysis cannot be reused automatically in these cases. In the degree three case the angular polynomials are $-\gamma u^3$ in the horizontal chart and $\gamma u$ in the vertical chart, adding a vertical equilibrium. The physical zero regularization limit also imposes a relation between $\gamma$ and $\beta$.

\subsection{Parameter diagram and finite-equilibrium regimes}\label{sec:paramdiagram}
\begin{figure}
\centering
\includegraphics[width=0.94\linewidth,height=0.74\textwidth,keepaspectratio]{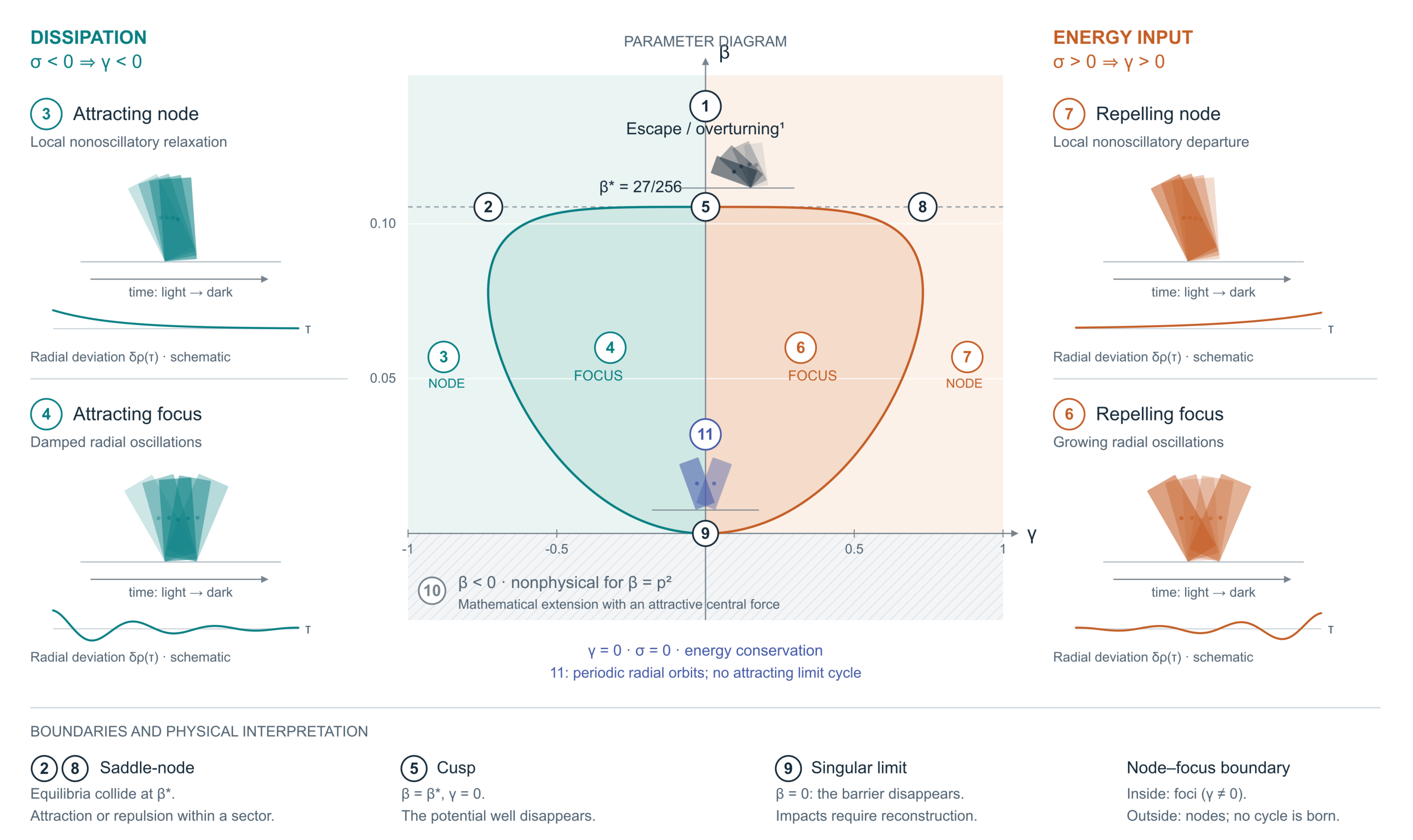}
\caption{Physical regime map of the reduced radial model, with case numbers 1--11 used throughout the text and Figure~\ref{fig:phasecomputed}. The curve is the finite equilibrium node--focus boundary \eqref{eq:nodefocus}. Block silhouettes illustrate radial trends, not reconstructed block trajectories.}
\label{fig:cardiofases1}
\end{figure}

\begin{figure}[p]
\centering
\includegraphics[width=\linewidth,height=0.83\textheight,keepaspectratio]{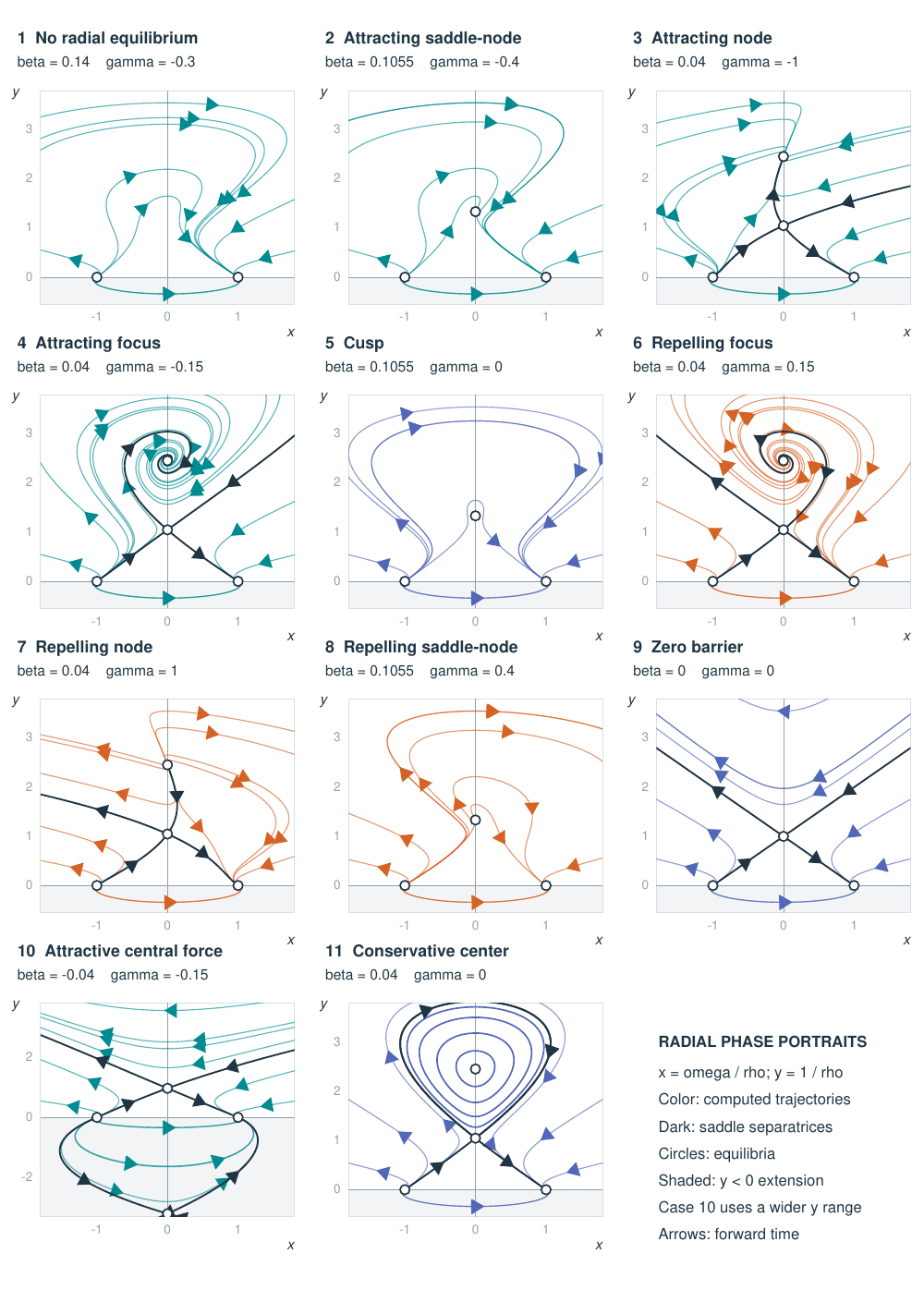}
\caption{Computed planar examples for cases 1--11 of Figure~\ref{fig:cardiofases1}, in coordinates $(x,y)=(\omega/\rho,1/\rho)$. Parameter values are printed above each panel. Arrows follow the polynomial field. Dark curves approximate saddle separatrices, and the closed center curves are energy levels. These are representative finite windows, not complete portraits at infinity.}
\label{fig:phasecomputed}
\end{figure}

Figure~\ref{fig:cardiofases1} summarizes the finite equilibrium regimes. For the equilibrium $y_e>4/3$, the node--focus boundary admits the parametrization
\begin{equation}\label{eq:nodefocus}
\beta=\frac{y_e-1}{y_e^4},\qquad
\gamma=\pm\frac{2\sqrt{3y_e-4}}{y_e^2},\qquad y_e>\frac43.
\end{equation}
Within the enclosed region the equilibrium is a focus for $\gamma\neq0$ and a center for $\gamma=0$. Outside it the equilibrium is a node, provided $0<\beta<\beta^*$. Finally, on the boundary it is an improper node. This change describes the local form of relaxation or repulsion, not the creation of a periodic attractor. Figure~\ref{fig:phasecomputed} provides computed planar examples of the cases listed at Figure~\ref{fig:cardiofases1}. The compactified radial schemes in Figures~\ref{fig:compactifieda} and~\ref{fig:compactifiedb} use the identifiers C1--C17, related to cases 1--11 by Table~\ref{tab:portraitcrosswalk}.

\subsection{Computed compactified radial portraits}\label{sec:compactified}
Figures~\ref{fig:compactifieda} and~\ref{fig:compactifiedb} show seventeen computed parameter examples. The panel identifiers C1--C17 refer to the parameters printed in each panel. They do not constitute an exhaustive topological classification. Each panel also gives its finite regime case number from Figures~\ref{fig:cardiofases1} and~\ref{fig:phasecomputed}, as summarized in Table~\ref{tab:portraitcrosswalk}.

The radial variables are projected to the Poincar\'e disk by
\begin{equation}
(\xi,\eta)=\frac{(\rho,\omega)}{\sqrt{1+\rho^2+\omega^2}}.
\end{equation}
Multiplication of the radial field by $\rho^4$, followed by degree five compactification, preserves orbit orientation on each open half disk. A further positive normalization is used for trajectory integration. The central diameter $\xi=0$ is excluded: the multiplier does not define physical impact continuation there. The outer circle represents radial infinity, rather than polynomial infinity in $(x,y)$ coordinates. Boundary arrows follow the compactified fiel. Additional equilibrium points introduced along the singular diameter by the multiplier should not be considered as physical equilibria.

The examples distinguish node and focus behavior and sample both signs of $\gamma^2-4\beta$. Conservative closed curves are obtained from the effective energy, and saddle separatrices are integrated from local eigendirections.

\begin{table}[ht]
\centering\small
\caption{Correspondence between the computed examples C1--C17 and the physical regime map. Parameters are printed in Figures~\ref{fig:compactifieda} and~\ref{fig:compactifiedb}. This table describes finite equilibrium regimes, not a complete global classification.}
\label{tab:portraitcrosswalk}
\begin{tabularx}{\linewidth}{@{}l l X@{}}
\toprule
Panel & Map case & Finite radial behavior \\
\midrule
C1--C3 & 1 & No finite radial equilibrium, $\beta>\beta^*$.\\
C4--C5 & 2 & Saddle--node with an attracting sector.\\
C6 & 5 & Conservative cusp.\\
C7--C8 & 8 & Saddle--node with a repelling sector.\\
C9 & 3 & Attracting node and saddle.\\
C10, C12 & 4 & Attracting focus and saddle.\\
C11 & 3 & Attracting node with $\gamma^2<4\beta$.\\
C13 & 11 & Center and saddle; conservative periodic family.\\
C14, C16 & 6 & Repelling focus and saddle.\\
C15 & 7 & Repelling node with $\gamma^2<4\beta$.\\
C17 & 7 & Repelling node and saddle.\\
Not shown & 9, 10 & $\beta=0$ and $\beta<0$ are not covered by the 17 computed examples.\\
\bottomrule
\end{tabularx}
\end{table}

\begin{figure}[p]
\centering
\includegraphics[width=\linewidth]{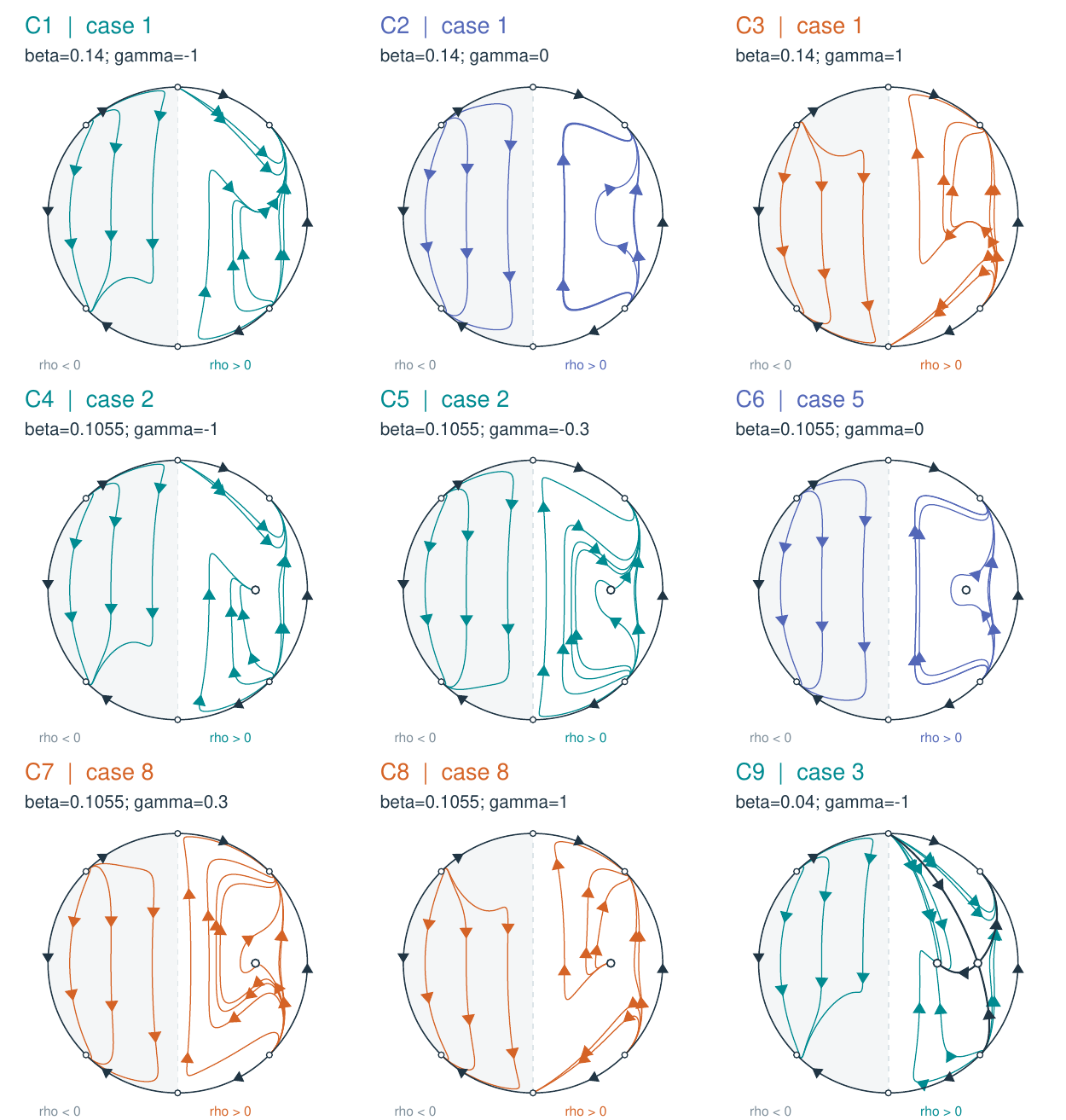}
\renewcommand{\thefigure}{\arabic{figure}a}
\caption{Computed compactified radial portraits, C1--C9. Each panel prints its parameters and corresponding finite regime case. Colored curves are radial trajectories. Dark curves are saddle separatrices, when present. The right half is $\rho>0$, the dashed diameter is the excluded singular set $\rho=0$, and the outer circle is radial infinity. Arrows preserve physical time orientation. Panels C10--C17 appear in Figure~\ref{fig:compactifiedb}.}
\label{fig:compactifieda}
\end{figure}
\begin{figure}[p]
\centering
\includegraphics[width=\linewidth]{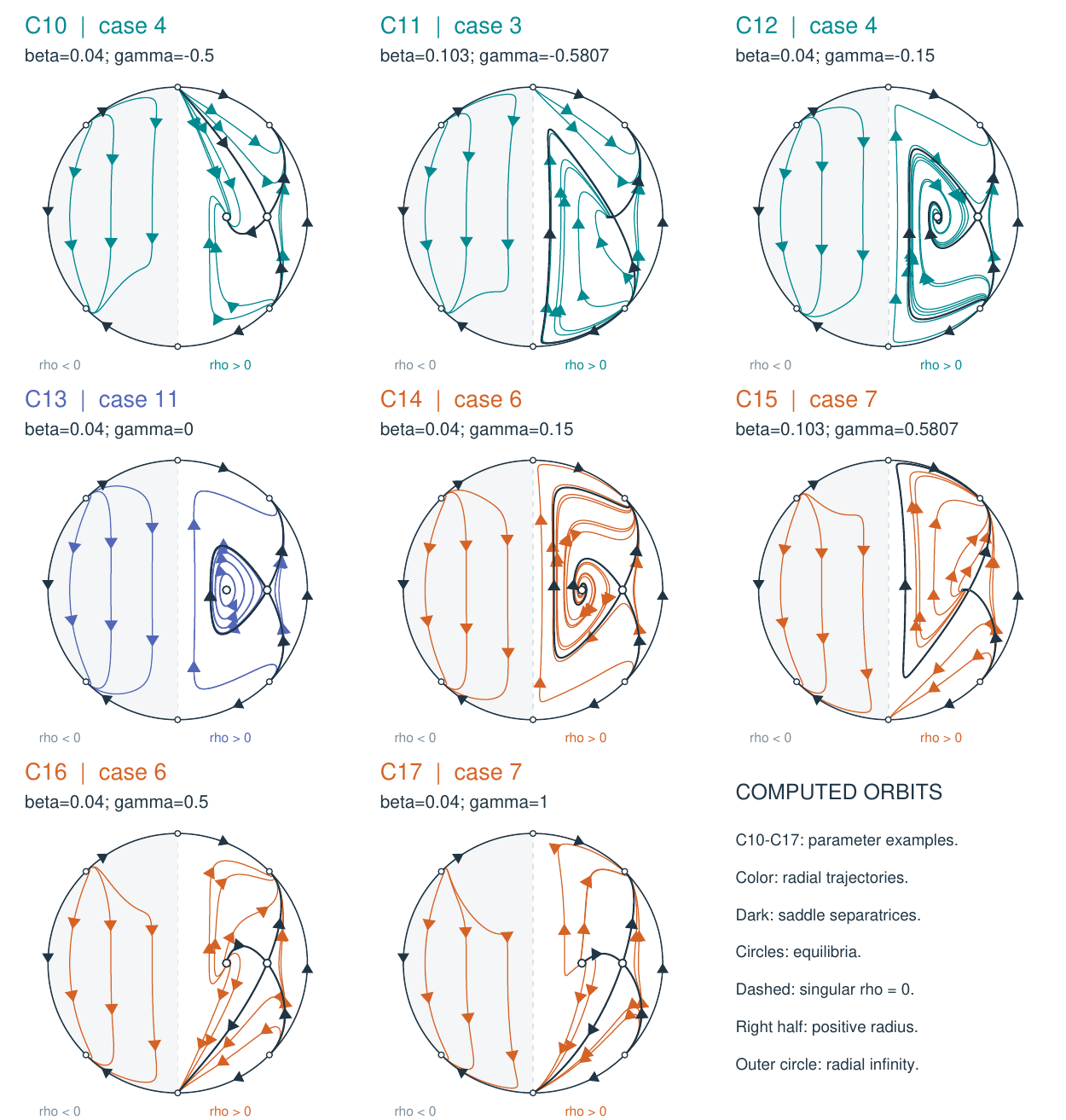}
\addtocounter{figure}{-1}
\renewcommand{\thefigure}{\arabic{figure}b}
\caption{Computed compactified radial portraits, C10--C17, with the conventions of Figure~\ref{fig:compactifieda}. The conservative closed curves in C13 are energy levels. C11 and C15 are explicit node examples with $\gamma^2<4\beta$ (some of their equilibria are shown with reduced size to improve readability). These finite numerical examples are not an exhaustive bifurcation classification at the singular boundary.}
\label{fig:compactifiedb}
\end{figure}

\clearpage
\section{Physical interpretation of the mathematical results}\label{sec:physical}
\subsection{What the phase-plane coordinates measure}
Radial variable $r=\alpha\rho$ is the modulus of complexified rocking coordinate. From \eqref{eq:physicalcoordinates}, we have that $r=\alpha/y$ and $x=\rho'/\rho$. Hence, $x>0$ corresponds to radial expansion while $x<0$ corresponds to radial contraction; they do not describe rocking to the right or left. Physical radial range of interest is $y>0$, while the lower half-plane is an artificial extension. The nodes $P_1=(-1,0)$ and $P_2=(1,0)$ are situated at infinite radius and not at finite block configurations. Attraction to $P_2$ thus symbolizes radial escape and not mechanical stabilization.

The original parametrization (see Table~\ref{tab:symbols} for reference) becomes
\begin{equation}\label{eq:physicalparameters}
\beta=p^2,\qquad \gamma=\sigma\sqrt\beta,\qquad
\sigma=\frac{2}{\pi}\ln\mu.
\end{equation}
Hence $0<\mu<1$ chooses $\gamma<0$, whereas $\mu=1$ -no damping- renders $\gamma=0$. A positive value for $\gamma$ means an addition of energy at impact (as in the jumping robot case but not in the block example). On the other hand, $\beta<0$ has no physical meaning. For a constant restitution coefficient the restriction $\gamma=\sigma\sqrt\beta$ in the two-parameter graph, bounds $\gamma$ and $\beta$ (Fig.~\ref{fig:cardiofases1}). When $\mu=0$ -infinite damping-, $\gamma\to -\infty$.
\subsection{An effective well, an escape barrier, and an energy balance}
From the dimensionless effective radial potential and energy by
\begin{equation}\label{eq:effectiveenergy}
V_\beta(\rho)=-\frac{\rho^2}{2}+\rho+\frac{\beta}{2\rho^2},\qquad
E=\frac{\omega^2}{2}+V_\beta(\rho),
\end{equation}
the radial equation is $\rho''=-V_\beta'(\rho)+\gamma\rho'/\rho^2$. Therefore
\begin{equation}\label{eq:energybalance}
\frac{dE}{d\tau}=\gamma\frac{\omega^2}{\rho^2}.
\end{equation}
This is the dissipated power of the reduced small-angle model, up to an additive constant and the dimensional scale $MgR\alpha^2$. The $\beta/(2\rho^2)$ term is the centrifugal barrier generated by conservation of the complex phase momentum. The gravitational part opposes small radial excursions but favors outward motion beyond its barrier.

Stationary radii satisfy
\begin{equation}\label{eq:radialequilibria}
\beta=\rho_e^3(1-\rho_e),\qquad
V_\beta''(\rho_e)=\frac{3}{\rho_e}-4.
\end{equation}
There exist two values of $\rho_e$, a minimum value $\rho_m \in (0,3/4)$ and a maximum value $\rho_s \in (3/4,1)$ for $0<\beta < 27/256$. The minimum value is the node, focus or center and the maximum value is the saddle. Their collision at $\rho_*=3/4$ is what gives rise to the critical value of $\beta = \beta^* = 27/256$ because the well and the escape barrier vanish simultaneously. For $\beta > \beta^*$, the force is always positive.
For $\gamma<0$, energy decreases whenever $\omega\neq0$. A useful sufficient condition for capture is
\begin{equation}\label{eq:capture}
0<\rho_0<\rho_s,\qquad E(\rho_0,\omega_0)<E_s:=V_\beta(\rho_s).
\end{equation}
The trajectory cannot reach $\rho_s$ because that would require energy at least $E_s$, nor can it reach $\rho=0$ because $V_\beta\to+\infty$ there. It remains in a compact region and converges to $(\rho_m,0)$: the only invariant set there on which $E'=0$ is that equilibrium. This also rules out a finite-energy collision with $\rho=0$ for $\beta>0$, $\gamma\leq0$. Condition~\eqref{eq:capture} is sufficient, not necessary; dissipation can capture some trajectories initially above the barrier. The saddle's stable separatrices provide the dynamical capture/escape threshold, which is not a constant-energy curve when $\gamma\neq0$.

For $\gamma=0$, energies between $V_\beta(\rho_m)$ and $E_s$ have a bounded periodic component within the well. At $E_s$, its boundary is the homoclinic loop to the saddle. In case that $\gamma\neq0$, system~\eqref{sys:Rockinggeneral} does not have any nontrivial radial periodic solution. Indeed, since over a period $T$ the total energy returns to its initial state $E(T)=E(0)$, while integrating \eqref{eq:energybalance}, we obtain 
\begin{equation*}
E(T)-E(0)=\gamma\int_0^T\frac{\omega^2}{\rho^2}\,d\tau\neq0,
\end{equation*}
since the integrand is nonnegative and vanishes identically only on an equilibrium. In particular, positive energy input at impact ($\gamma>0$) cannot by itself sustain a stable oscillation of fixed amplitude: a limit cycle would require gain and loss to balance along the orbit, but $dE/d\tau$ never changes sign. Trajectories instead spiral away from the minimum.
\subsection{Physical reading of the eleven portraits in the $(\gamma,\beta)$ plane}\label{sec:portraits}
\begin{description}
\item[Case 1.] For $\beta>\beta^*$ there is no potential well. Outgoing trajectories are unbounded; the phase plane has no finite fixed point that supports radial dynamics.
\item[Cases 2 and 8.] At $\beta=\beta^*$ and $\gamma\neq0$, saddle-node exhibit an attracting sector if $\gamma<0$ and repelling if $\gamma>0$. It is a breakdown of the well, but not a stable structure under any perturbations.
\item[Cases 3 and 4.] If $\gamma<0$, the minimum attracts the domain as a node or focus. In the former case, we have non-oscillatory radial relaxation; in the latter, oscillatory one with damping. Neither of the dynamics guarantees capture of the domain.
\item[Case 5.] When $\beta=\beta^*$ and $\gamma=0$, the center and the saddle merge at the nilpotent cusp and destroy the conservative period annulus.
\item[Cases 6 and 7.] When $\gamma >0$, the minimum is the repelling focus or node. In relation to that, radial perturbations increase exponentially with or without local spirals.
\item[Case 9.] When $\beta=0$, the existing finite radial equilibrium becomes the saddle $\rho=1$. Without the centrifugal term, no impact continuation law exists at $\rho=0$.
\item[Case 10.] For $\beta<0$, the polynomial plane has two saddle points, only one for $y>0$. The inverse cubic force term that attracts is part of the extension and not the original momentum-block model.
\item[Case 11.] For $\gamma = 0$ and $0 < \beta < \beta^*$, a periodic orbit surrounds the center up to the homoclinic point. The energy level and the amplitude of motion depend on the initial conditions and none of these attract their neighbors.
\end{description}

Linearization around the minimum makes the difference between nodes and foci more explicit. Taking $\eta=\rho-\rho_m$:
\begin{equation}\label{eq:radiallinearization}
\eta''+c_m\eta'+\Omega_m^2\eta=0,\qquad
c_m=-\frac{\gamma}{\rho_m^2},\qquad
\Omega_m^2=\frac{3}{\rho_m}-4.
\end{equation}
For $\gamma<0$, we have $c_m>0$: the cases 3 and 4 represent the overdamped and underdamped radial motions, with critical damping in \eqref{eq:nodefocus}. Here we obtain $\tau$-dependent frequencies and decay rates; their dimensional counterparts must be scaled by $\sqrt{MgR/I}$.
\subsection{Radial stabilization is a relative equilibrium}
For any $p\neq0$, a radial steady state satisfies
\begin{equation}
\rho=\rho_m,\qquad \omega=0,\qquad \psi'=\frac{p}{\rho_m^2}.
\end{equation}
The rotating complex variable has a fixed modulus. Thus, the radially attracting equilibrium is a relative equilibrium for the extended system $(\rho,\psi)$ and not a block at rest. Indeed, at this point the damping force vanishes since it acts only through radial velocity. The phase momentum $p$ is conserved, which justifies the persistence of non-zero phase motion in spite of $\gamma<0$. Likewise, a periodic radial orbit of period $T$ produces the phase advance
\begin{equation}
\Delta\psi=p\int_0^T\frac{d\tau}{\rho(\tau)^2}.
\end{equation}
Closure of the entire motion is only possible when $\Delta\psi/(2\pi)$ is rational. Else, closure of the radial period motion occurs but there is no closure of the complex motion. The reconstruction of signed rocking motion requires both radial solutions and phase information, together with the impact limit criterion.
\section{Conclusions}
Regarding the polynomial representation, it is now possible to establish the direct relationship between the geometric properties of phase portraits and the mechanical properties of regularized radial dynamics. The case of positive radius can be considered as a consequence of the transformation $x=\omega/\rho$, $y=1/\rho$. It keeps the same dimensionless time and yields the polynomial system for the singular radial dynamics. Within this kind of equivalence, it is now possible to apply the planar analysis of finite equilibrium points and separatrices. Now the meaning of the new boundaries becomes obvious. For instance, the attraction to $y=0$ implies radial escaping rather than stabilization at the finite blocks.\\
Moreover, the effective potential shows what makes confinement possible. For $0<\beta<27/256$, centrifugal regularization gives rise to a well and an inner barrier, whereas gravity contributes an outer saddle barrier, and their intersection at $\rho=3/4$ annihilates the well when $\beta^*=\beta=27/256$. Such a value of $\beta^*$ marks the limit of existence of confined radial states but does not necessarily dictate the final result of any particular set of initial conditions. In the case of dissipation, initial states lying within the capture domain inside the saddle and having the energy smaller than the barrier height will be captured and flow into the minimum point. Outside the above domain, stable separatrices of the saddle, and not just the type of equilibrium itself, make the difference between capture and escape.

The equation $E'=\gamma\omega^2/\rho^2$ defines the energetic difference between the two main regimes. Dissipation creates radial relaxation, conservation allows for a range of periodic solutions inside the potential well, and energy input leads to instability of the potential minimum. The node-focus separatrix defines the transition from nonoscillatory and oscillatory radial behavior, respectively. Note that the node-focus dynamics are not an oscillatory solution by themselves. If $\gamma\neq0$, monotonic energy means there are no periodic radial solutions at all. For $\gamma=0$, periodic radial solutions do exist, but they are not isolated: each energy level between $V_\beta(\rho_m)$ and $E_s$ is a closed orbit, so they form a continuous family (the period annulus) rather than a limit cycle, and their amplitude is fixed by the initial conditions.\\

We stress that the phase must be reconstructed before these results are translated to the rocking block. At radial equilibrium with finite $p\neq 0$, $\rho=\rho_m$ and $\omega=0$, but the phase keeps rotating, $\psi'=p/\rho_m^2\neq0$, because $p$ is conserved and the damping acts only through the radial velocity. Now variable $\theta=\rho_m e^{i\psi}$ moves uniformly on a circle. This is a relative equilibrium of the regularized system in the complex $(\rho,\psi)$ plane, not a block at rest; the physical coordinate $\operatorname{Re}\theta=\rho_m\cos\psi$ oscillates with amplitude $\rho_m$, which vanishes only in the limit $p\to0$. Also, along a periodic radial orbit of period $T$ the phase advances by $\Delta\psi=p\int_0^T\rho^{-2}\,d\tau$, and the complex motion is periodic only if $\Delta\psi/(2\pi)$ is rational; otherwise it would be quasi-periodic. Radial attraction, periodic radial motion and a block at rest are different scenarios. 

This finite-to-equilibrium transition should also be kept distinct from the boundary discriminant $\gamma^2=4\beta$. This discriminant applies only to singular sectors uncovered via compactification and blowup, not to damping towards the potential minimum. On a family of fixed restitution, $\gamma=\sigma\sqrt\beta$; hence the physical parameters trace a constrained trajectory in the two-parameter diagram. For a fixed value of $\beta>0$ and $\gamma\leq0$, the energy barrier continues to preclude a finite-energy collision with $\rho=0$, regardless of the value of this discriminant. Boundary lines that are compactified and have an unbounded velocity cannot be viewed as finite-velocity collisions.\\

Finally, notice that radial escape indicates loss of confinement in the small angle model only (quantitative overturning for larger angles and times would require additional analysis that we omit here). As in \cite{Meccanica}, to recover \textit{pure} physical rocking, one requires the limit $(\beta,\gamma)=(p^2,p\sigma)\to(0,0)$ along with phase jump reconstruction. However, setting the parameters strictely to zero  in the complex regularization would destroy the concentrated dissipation near the impact. In this regard, the model captures increasingly more accurate concentrated impacts when $p\to 0$.
\clearpage
\begin{table}[!ht]
\centering\small
\caption{Physical and reduced-model notation.}
\label{tab:symbols}
\begin{tabularx}{\textwidth}{@{}lX@{}}
\toprule
Symbol & Meaning \\
\midrule
$I,M,R,g$ & Pivot inertia, mass, distance from pivot to center of mass, and gravitational acceleration.\\
$\alpha,\theta$ & Slenderness angle and signed physical rocking angle.\\
$r,\psi$ & Positive modulus and phase of the complexified coordinate.\\
$\mu,\sigma$ & Restitution coefficient and coupling $\sigma=(2/\pi)\ln\mu$, for $\mu>0$.\\
$k,\tau$ & $k=MgR/I$ and dimensionless time $\tau=t\sqrt{k}$.\\
$p_\psi,p$ & Conserved phase momentum $Ir^2\dot\psi$ and its dimensionless value $p=p_\psi/(\alpha^2\sqrt{MgRI})$.\\
$\rho,\omega$ & Normalized radius $r/\alpha$ and radial velocity $d\rho/d\tau$.\\
$\beta,\gamma$ & Generalized parameters; the physical family has $\beta=p^2$, $\gamma=p\sigma$ with $p>0$.\\
$\beta^*$ & Well-disappearance threshold $27/256$.\\
$V_\beta,E$ & Effective potential and radial energy in \eqref{eq:effectiveenergy}.\\
$\rho_m,\rho_s$ & Minimum and maximum of $V_\beta$ for $0<\beta<\beta^*$.\\
$E_s$ & Saddle-barrier energy $V_\beta(\rho_s)$.\\
\bottomrule
\end{tabularx}
\end{table}

\begin{table}[!ht]
\centering\small
\caption{Polynomial and compactification notation.}
\begin{tabularx}{\textwidth}{@{}lX@{}}
\toprule
Symbol & Meaning \\
\midrule
$x,y$ & Polynomial coordinates $x=\omega/\rho$, $y=1/\rho$.\\
$P,Q$ & Right-hand sides of \eqref{sys:nuestracartaU1}.\\
$P_1,P_2$ & Boundary nodes $(-1,0)$ and $(1,0)$, representing infinite radius.\\
$P_3,P_4$ & Finite equilibria on $x=0$, when present.\\
$u,v$ & Compactification coordinates $(x,y)=(1/v,u/v)$.\\
$w,z$ & Weighted blow-up variables $u=w$, $v=w^2z$.\\
$\Delta$ & Exceptional-divisor discriminant $\gamma^2-4\beta$; distinct from the finite-equilibrium discriminant.\\
$\mathcal H$ & Conservative first integral, equal to $2E$ in polynomial coordinates.\\
$\Theta,\delta$ & Heaviside function and Dirac distribution.\\
\bottomrule
\end{tabularx}
\end{table}

\section*{Declaration of competing interest}

\section*{Acknowlwdgements}
 The authors are supported by the Agencia Estatal de Investigación of Spain grant PID2023-147461NB-I00.

\section*{Data availability}
No data was used for the research described in the article.
\clearpage
\bibliography{rm}

@article{Dirac,
   author = {F. Prieto and P. B. Louren{\c{c}}o and C. S. Oliveira},
   doi = {10.1002/eqe.381},
   issn = {00988847},
   issue = {7},
   journal = {Earthquake Engineering and Structural Dynamics},
   pages = {839-857},
   publisher = {John Wiley and Sons Ltd},
   title = {Impulsive Dirac-delta forces in the rocking motion},
   volume = {33},
   year = {2004},
}

@article{Meccanica,
 author = {F. Prieto and P. B. Louren{\c{c}}o},
   doi = {10.1007/s11012-004-5875-z},
   journal = {Meccanica},
   pages = {121-133},
   title = {On the Rocking Behavior of Rigid Objects},
   year = {2005},
}

@article{Pena,
   author = {F. Pe{\~n}a and F. Prieto and P. B. Louren{\c{c}}o and A. Campos Costa and J. V. Lemos},
   doi = {10.1002/eqe.739},
   issn = {10969845},
   issue = {15},
   journal = {Earthquake Engineering and Structural Dynamics},
   pages = {2383-2399},
   publisher = {John Wiley and Sons Ltd},
   title = {On the dynamics of rocking motion of single rigid-block structures},
   volume = {36},
   year = {2007},
}

@article{Chattering2001,
   author = {D J Wagg and S R Bishop},
   issue = {1},
   journal = {International Journal of Bifurcation and Chaos},
   pages = {57-71},
   title = {CHATTER, STICKING AND CHAOTIC IMPACTING MOTION IN A TWO-DEGREE OF FREEDOM IMPACT OSCILLATOR},
   volume = {11},
   url = {www.worldscientific.com},
   year = {2001},
}

@article{Hogan,
   author = {S. J. Hogan},
   doi = {10.1098/rspa.1989.0114},
   issn = {0080-4630},
   issue = {1869},
   journal = {Proceedings of the Royal Society of London. A. Mathematical and Physical Sciences},
   month = {10},
   pages = {441-476},
   publisher = {The Royal Society},
   title = {On the dynamics of rigid-block motion under harmonic forcing},
   volume = {425},
   year = {1989},
}

@inproceedings{Makris2013,
   author = {Nicos Makris and Michalis F. Vassiliou},
   doi = {10.7712/120113.4702.c1083},
   journal = {ECCOMAS Thematic Conference - COMPDYN 2013: 4th International Conference on Computational Methods in Structural Dynamics and Earthquake Engineering, Proceedings - An IACM Special Interest Conference},
   pages = {2764-2785},
   publisher = {National Technical University of Athens},
   title = {Rocking response and stability analysis of an array of free-standing columns capped with a free-standing rigid beam},
   year = {2013},
}

@article{Proteins2015,
   author = {Peixiang Ma and Yi Xue and Nicolas Coquelle and Jens D. Haller and Tairan Yuwen and Isabel Ayala and Oleg Mikhailovskii and Dieter Willbold and Jacques Philippe Colletier and Nikolai R. Skrynnikov and Paul Schanda},
   doi = {10.1038/ncomms9361},
   issn = {20411723},
   journal = {Nature Communications},
   month = {10},
   pmid = {26436197},
   publisher = {Nature Publishing Group},
   title = {Observing the overall rocking motion of a protein in a crystal},
   volume = {6},
   year = {2015},
}

@article{Proteins2023,
   author = {Alexey Krushelnitsky and G{\"u}nter Hempel and Hannes Jurack and Tiago Mendes Ferreira},
   doi = {10.1039/d3cp00444a},
   issn = {14639076},
   issue = {23},
   journal = {Physical Chemistry Chemical Physics},
   month = {6},
   pages = {15885-15896},
   pmid = {37259579},
   publisher = {Royal Society of Chemistry},
   title = {Rocking motion in solid proteins studied by the 15N proton-decoupled {R1$\rho$} relaxometry},
   volume = {25},
   year = {2023},
}

@article{MetaInterfaces,
   author = {N. Hima and F. D'Annibale and F. Dal Corso},
   doi = {10.1016/j.ijmecsci.2022.108005},
   issn = {00207403},
   journal = {International Journal of Mechanical Sciences},
   month = {3},
   publisher = {Elsevier Ltd},
   title = {Non-smooth dynamics of buckling based metainterfaces: Rocking-like motion and bifurcations},
   volume = {242},
   year = {2023},
}

@inproceedings{Walking2019,
   author = {Qian Xiang and Kazunobu Hashimoto and Lei Sun and Jianjun Zhang and Shijie Guo},
   doi = {10.1109/ICARM.2019.8834245},
   isbn = {9781728100647},
   journal = {2019 4th IEEE International Conference on Advanced Robotics and Mechatronics,ICARM 2019},
   month = {7},
   pages = {862-867},
   publisher = {Institute of Electrical and Electronics Engineers Inc.},
   title = {An improved rocker-based inverted pendulum walk model for human gait},
   year = {2019},
}

@inproceedings{Walking2006,
   author = {J. M. Bourgeot and C. Canudas-De-Wit and B. Brogliato},
   doi = {10.1007/3-540-26415-9-61},
   isbn = {3540264132},
   journal = {Proceedings of the 8th International Conference on Climbing and Walking Robots and the Support Technologies for Mobile Machines, CLAWAR 2005},
   pages = {509-516},
   publisher = {Kluwer Academic Publishers},
   title = {Impact shaping for double support walk: From the rocking block to the biped Robot},
   year = {2006},
}

@article{AI,
   author = {Seyed Amir Banimahd and Anastasios I Giouvanidis and Shaghayegh Karimzadeh and Paulo B Louren{\c{c}}o},
   doi = {10.1002/eqe.4110},
   journal = {Earthquake Engineering and Structural Dynamics},
   title = {A multi-level approach to predict the seismic response of rigid rocking structures using artificial neural networks},
   url = {https://onlinelibrary.wiley.com/doi/10.1002/eqe.4110},
   year = {2024},
}

@article{Rega2006,
 author = {Stefano Lenci and Giuseppe Rega},
   doi = {10.1016/J.CHAOS.2005.07.007},
   issn = {0960-0779},
   issue = {2},
   journal = {Chaos, Solitons and Fractals},
   month = {4},
   pages = {527-542},
   publisher = {Pergamon},
   title = {A dynamical systems approach to the overturning of rocking blocks},
   volume = {28},
   year = {2006},
}

@article{MolecularRotor,
   author = {Cosima Stahler and Daisy R.S. Pooler and Romain Costil and Dhruv Sudan and Pieter van der Meulen and Ryojun Toyoda and Ben L. Feringa},
   doi = {10.1021/ACS.JOC.2C01830/SUPPL-FILE/JO2C01830-SI-002.ZIP},
   issn = {15206904},
   issue = {1},
   journal = {Journal of Organic Chemistry},
   month = {1},
   pages = {1-8},
   publisher = {American Chemical Society},
   title = {Coupled Rocking Motion in a Light-Driven Rotary Molecular Motor},
   volume = {89},
   url = {https://pubs.acs.org/doi/full/10.1021/acs.joc.2c01830},
   year = {2024},
}

@book{Raibert1986,
  author = {Raibert, Marc H.},
  title = {Legged Robots That Balance},
  publisher = {MIT Press},
  year = {1986},
  address = {Cambridge, MA}
}

@article{GarciaPerezChavelaSusin2006,
  author = {Garc{\'i}a, A. and P{\'e}rez-Chavela, E. and Susin, A.},
  title = {A generalization of the {P}oincar{\'e} compactification},
  journal = {Archive for Rational Mechanics and Analysis},
  volume = {179},
  pages = {285--302},
  year = {2006}
}

@book{DumortierLlibreArtes2006,
  author = {Dumortier, F. and Llibre, J. and Art{\'e}s, J. C.},
  title = {Qualitative Theory of Planar Differential Systems},
  series = {Universitext},
  publisher = {Springer-Verlag},
  address = {Berlin},
  year = {2006}
}

\end{document}